\documentclass[12pt,twoside]{article}
\usepackage{amsmath, amssymb, amsthm,graphicx}
\font\fourteenb=cmb10 at 14pt
\begin{document}

\vspace*{-1.0cm}\noindent \copyright
 Journal of Technical University at Plovdiv\\[-0.0mm]\
\ Fundamental Sciences and Applications, Vol. 2, 1996\\[-0.0mm]
\textit{Series A-Pure and Applied Mathematics}\\[-0.0mm]
\ Bulgaria, ISSN 1310-8271\\[+1.2cm]
\font\fourteenb=cmb10 at 12pt
\begin{center}

{\bf \LARGE Note on the multipliers of Cauchy integrals of
logarithmic potentials
   \\ \ \\ \large Peyo Stoilov}
\end{center}

\footnotetext{{\it Received June 15,1996.}}
\begin{abstract}
The present note contains a generalization of a theorem of
Hallenbeck and Samotij for the multipliers of Cauchy integrals of
logarithmic potentials.
\end{abstract}

\section{Introduction.}

     Let  $D$  denote the unit disk in the complex plane and $T$  the
unit circle. Let  $M$  be the Banach space of all complex-valued
Borel  measures on  $T$  with the usual variation norm. For
$\alpha \ge 0$ ,  let  $F_{\alpha } $  denote the family of the
analytic functions  $g$ ,  for which there exists  $\mu \in M$
such that

$$\displaystyle g(z)=\int _{T}\frac{1}{(1-\overline{\varsigma
}z)^{\alpha } }  {\kern 1pt} {\kern 1pt} {\kern 1pt} {\kern 1pt}
{\kern 1pt} d\mu (\varsigma )=K_{\alpha }^{\mu } ,{\kern 1pt}
{\kern 1pt} {\kern 1pt} {\kern 1pt} {\kern 1pt} {\kern 1pt} {\kern
1pt} {\kern 1pt} {\kern 1pt} {\kern 1pt} {\kern 1pt} {\kern 1pt}
{\kern 1pt} {\kern 1pt} {\kern 1pt} {\kern 1pt} {\kern 1pt} {\kern
1pt} {\kern 1pt} {\kern 1pt} {\kern 1pt} {\kern 1pt} {\kern 1pt}
{\kern 1pt} {\kern 1pt} {\kern 1pt} {\kern 1pt} {\kern 1pt} {\kern
1pt} {\kern 1pt} {\kern 1pt} {\kern 1pt} {\kern 1pt} {\kern 1pt}
{\kern 1pt} {\kern 1pt} {\kern 1pt} {\kern 1pt} {\kern 1pt} {\kern
1pt} {\kern 1pt} {\kern 1pt} {\kern 1pt} {\kern 1pt} {\kern 1pt}
{\kern 1pt} {\kern 1pt} {\kern 1pt} {\kern 1pt} {\kern 1pt} {\kern
1pt} {\kern 1pt} {\kern 1pt} {\kern 1pt} {\kern 1pt} {\kern 1pt}
{\kern 1pt} {\kern 1pt} {\kern 1pt} {\kern 1pt} {\kern 1pt} {\kern
1pt} {\kern 1pt} {\kern 1pt} {\kern 1pt} {\kern 1pt} {\kern 1pt}
{\kern 1pt} {\kern 1pt} {\kern 1pt} {\kern 1pt} {\kern 1pt} {\kern
1pt} {\kern 1pt} {\kern 1pt} {\kern 1pt} {\kern 1pt} {\kern 1pt}
\alpha >0,$$

$$\displaystyle g(z)=g(0)+\int _{T}\log
\left(\frac{1}{1-\overline{\varsigma }z} \right) {\kern 1pt}
{\kern 1pt} {\kern 1pt} {\kern 1pt} {\kern 1pt} d\mu (\varsigma
)=K_{0}^{\mu } ,{\kern 1pt} {\kern 1pt} {\kern 1pt} {\kern 1pt}
{\kern 1pt} {\kern 1pt} {\kern 1pt} {\kern 1pt} {\kern 1pt} {\kern
1pt} {\kern 1pt} {\kern 1pt} {\kern 1pt} {\kern 1pt} {\kern 1pt}
{\kern 1pt} {\kern 1pt} {\kern 1pt} {\kern 1pt} {\kern 1pt} {\kern
1pt} {\kern 1pt} {\kern 1pt} {\kern 1pt} {\kern 1pt} {\kern 1pt}
\alpha =0.$$

     We  note  that   $F_{\alpha } $  is  a Banach  space with the natural norm

$$\displaystyle \left\| g\right\| _{F_{\alpha } } =\inf
\left\{\left\| \mu \right\| :\mu \in M,{\kern 1pt} {\kern 1pt}
{\kern 1pt} {\kern 1pt} {\kern 1pt} {\kern 1pt} {\kern 1pt} {\kern
1pt} {\kern 1pt} {\kern 1pt} {\kern 1pt} {\kern 1pt} g=K_{\alpha
}^{\mu } \right\}.$$

     Let  $m_{\alpha } $ denote the set of all multipliers of   $F_{\alpha } $  and

$$\displaystyle \left\| f\right\| _{m_{\alpha } } =\sup
\left\{\left\| f{\kern 1pt} g\right\| _{F_{\alpha } } :\left\|
g\right\| _{F_{\alpha } } \le 1\right\}.$$

\

\
   The following results were proved in [1].

\

\
     {\bf Theorem A}.   $m_{0} \subset m_{\alpha } \subset H^{\infty } $      $(m_{0} \ne m_{\alpha } )$ {\it     for}{\it  }{\it  each}   $\alpha >0$ .

\

\ {\bf Theorem B.}  {\it If } $f$ {\it is an analy}{\it tic
function on } $D$ {\it  }{\it and }{\it  } $f'\in H^{p} $ {\it
for some } $p>1$ {\it , }{\it  then } $f\in m_{0} $ {\it  .}

\

\ {\bf Theorem C.} {\it For every } $p\ge 1$ {\it  there{\kern
1pt}{\kern 1pt}}{\kern 1pt}{\it }{\it {\kern 1pt}{\kern
1pt}is}{\kern 1pt}{\it  }{\it {\kern 1pt}a{\kern 1pt}}{\kern
1pt}{\kern 1pt}{\it }{\it {\kern 1pt}constant }{\it  } $c_{p} $
{\it , such that }{\it }{\it for every } $f\in m_{0} $ {\it  and
each  } $\varsigma \in T$ {\it we have}

     $$\int _{0}^{1}\left((1-r)\log \frac{1}{1-r} \right) ^{p-1} {\kern 1pt} {\kern 1pt} {\kern 1pt} \left|f'(r\varsigma )\right|^{p} {\kern 1pt} {\kern 1pt} dr{\kern 1pt} {\kern 1pt} {\kern 1pt} \le {\kern 1pt} {\kern 1pt} {\kern 1pt} {\kern 1pt} {\kern 1pt} c_{p} {\kern 1pt} {\kern 1pt} {\kern 1pt} \left\| f\right\| _{m_{0} } .$$

\

\
     The present note contains a generalization  of  {\bf Theorem C}.
\section{Main result}
      {\bf Theorem D.} {\it If } $f\in m_{0} $ {\it    then }{\it }{\it
for each }{\it  } $\varsigma \in T$ {\it  }{\it  and  } $0\le r<1$
{\it  we have}

\

\ {\it a) } $\left|f'(r\varsigma )\right|{\kern 1pt} {\kern 1pt}
{\kern 1pt} (1-r)\log \frac{1}{1-r} {\kern 1pt} {\kern 1pt} \le
{\kern 1pt} {\kern 1pt} {\kern 1pt} {\kern 1pt} {\kern 1pt} {\kern
1pt} {\kern 1pt} \left\| f\right\| _{m_{0} } +\left\| f\right\|
_{H^{\infty } } <\infty $ {\it   ;}

\

\
      {\it b)}  $\int _{0}^{1}\left|f'(r\varsigma )\right|{\kern 1pt}  {\kern 1pt} {\kern 1pt} {\kern 1pt} {\kern 1pt} {\kern 1pt} dr{\kern 1pt} {\kern 1pt} {\kern 1pt} {\kern 1pt} {\kern 1pt} {\kern 1pt} {\kern 1pt} {\kern 1pt} \le {\kern 1pt} {\kern 1pt} {\kern 1pt} {\kern 1pt} {\kern 1pt} {\kern 1pt} {\kern 1pt} {\kern 1pt} {\kern 1pt} {\kern 1pt} {\kern 1pt} \frac{\pi }{2} {\kern 1pt} {\kern 1pt} {\kern 1pt} \left\| f\right\| _{m_{1} } <{\kern 1pt} {\kern 1pt} {\kern 1pt} {\kern 1pt} {\kern 1pt} {\kern 1pt} \infty .$

\

\

     \emph{Proof.  a)} Let   $f\in m_{0} $ {\it ,   } $g_{0} =\log \frac{1}{1-z} {\kern 1pt} {\kern 1pt} {\kern 1pt} {\kern 1pt} {\kern 1pt} $ {\it  and   } $\xi =r\varsigma \in D.$

     By Closed Graph Theorem  the multiplication operator   $F_{0} \backepsilon g\to fg\in F_{0} $   is bounded and
$\displaystyle \left\| f{\kern 1pt} g\right\| _{F_{\alpha } } \le
\left\| f\right\| _{m_{\alpha } } \left\| g\right\| _{F_{\alpha }
} $ .

     Since  $g_{0} (z{\kern 1pt} {\kern 1pt} {\kern 1pt} {\kern 1pt} \overline{\xi {\kern 1pt} }{\kern 1pt} /\left|\xi \right|)\in F_{0} $ {\it    }and     $\left\| g_{0} (z{\kern 1pt} {\kern 1pt} {\kern 1pt} {\kern 1pt} \overline{\xi {\kern 1pt} }{\kern 1pt} /\left|\xi \right|)\right\| _{F_{0} } =1$    for each    $\xi \in D$ , then

$$\left\| f(z){\kern 1pt} {\kern 1pt} {\kern 1pt} g_{0} (z{\kern 1pt} {\kern 1pt} {\kern 1pt} {\kern 1pt} \overline{\xi {\kern 1pt} }{\kern 1pt} /\left|\xi \right|)\right\| _{F_{0} } \le \left\| f\right\| _{m_{0} } $$

and

$$\displaystyle f(z){\kern 1pt} {\kern 1pt} {\kern 1pt} \log
(1/(1-z{\kern 1pt} {\kern 1pt} {\kern 1pt} {\kern 1pt}
\overline{\xi {\kern 1pt} }{\kern 1pt} /\left|\xi \right|))=\int
_{T}\log \left(\frac{1}{1-\overline{\varsigma }z} \right) {\kern
1pt} {\kern 1pt} {\kern 1pt} {\kern 1pt} {\kern 1pt} d\mu _{\xi }
(\varsigma ),$$ {\it}

where        $\mu _{\xi } \in M,{\kern 1pt} {\kern 1pt} {\kern
1pt} {\kern 1pt} {\kern 1pt} {\kern 1pt} {\kern 1pt} {\kern 1pt}
{\kern 1pt} \left\| \mu _{\xi } \right\| \le \left\| f\right\|
_{m_{0} } $ .

     Then

$$f'(z){\kern 1pt} {\kern 1pt} {\kern 1pt} \log (1/(1-z{\kern 1pt} {\kern 1pt} {\kern 1pt} {\kern 1pt} \overline{\xi {\kern 1pt} }{\kern 1pt} /\left|\xi \right|))=\int _{T}\frac{\overline{\varsigma }}{1-\overline{\varsigma }z}  {\kern 1pt} {\kern 1pt} {\kern 1pt} d\mu _{\xi } (\varsigma )-f(z)(1/(1-z{\kern 1pt} {\kern 1pt} {\kern 1pt} {\kern 1pt} \overline{\xi {\kern 1pt} }{\kern 1pt} /\left|\xi \right|))\cdot (\overline{\xi {\kern 1pt} }{\kern 1pt} /\left|\xi \right|).$$

     This  implies that

$$\displaystyle \left|f'(z)\right|{\kern
1pt} \left|{\kern 1pt} \log (1/(1-z{\kern 1pt}{\kern 1pt}
\overline{\xi {\kern 1pt} }{\kern 1pt} /\left|\xi
\right|))\right|(1-\left|z\right|)\le \left\| f\right\| _{m_{0} }
+\left\| f\right\| _{H^{\infty }}{\kern 1pt}
(1-\left|z\right|){\kern 1pt}(1/\left|1-z{\kern 1pt}{\kern 1pt}
\overline{\xi {\kern 1pt} }{\kern 1pt}/\left|\xi \right|\right|.$$
{\it }

     Substituting   $z=\xi $ {\it   }we obtain \emph{a)}.

     \emph{The proof of b)} can be found in [3] or [2].

     {\bf Corollary 1.}     {\it If   } $f\in m_{0} $ {\it  }{\it , then }{\it }{\it
{\kern 1pt}for{\kern 1pt}{\kern 1pt}each }{\it  } $\varsigma \in
T$ {\it }{\it and }{\it  } $p\ge 1$ {\it  }{\it  }{\it we have}

$$\displaystyle \int _{0}^{1}\left((1-r)\log \frac{1}{1-r} \right)
^{p-1} {\kern 1pt} {\kern 1pt} {\kern 1pt} \left|f'(r\varsigma
)\right|^{p} {\kern 1pt} {\kern 1pt} dr{\kern 1pt} {\kern 1pt}
{\kern 1pt} \le {\kern 1pt} {\kern 1pt} {\kern 1pt} {\kern 1pt}
{\kern 1pt} \frac{\pi }{2} {\kern 1pt} {\kern 1pt} {\kern 1pt}
\left\| f\right\| _{m_{1} } ({\kern 1pt} {\kern 1pt} {\kern 1pt}
\left\| f\right\| _{m_{0} } +\left\| f\right\| _{H^{\infty } }
)^{p-1} .$$

     The Proof it follows at once from \emph{a) and b)}.

     Let  $m(B)$ denote the family of all multipliers of the Bloch space B.

\

\
      {\bf Corollary 2. } $m_{0} \subset \left(m_{\alpha } \cap m(B)\right)$  {\it for each } $\alpha >0$ {\it .}

\

\

     \emph{Proof.} It follows at once from Theorem A and the assertion that  $f\in m_{0} $ {\it   }if  and only if  $f\in H^{\infty } $  and satisfies the condition \emph{a)} [4].

\noindent
{\small Department of Mathematics\\
        Technical University\\
        25, Tsanko Dijstabanov,\\
        Plovdiv, Bulgaria\\
        e-mail: peyyyo@mail.bg}

\end{document}